\documentclass[twoside,12pt,leqno]{article}
\usepackage{amsmath, amsthm, amsfonts, amssymb, color}
\usepackage{mathrsfs}
\usepackage{fancyhdr}
\usepackage{CJK}
\usepackage{bbm}
\usepackage{enumerate}
\usepackage{graphicx}
\usepackage{indentfirst}
\usepackage{stmaryrd}
\usepackage{refstyle}

\newcommand{\be}{\begin{eqnarray}}
\newcommand{\ee}{\end{eqnarray}}
\newcommand{\ce}{\begin{eqnarray*}}
\newcommand{\de}{\end{eqnarray*}}

\newtheorem{thm}{Theorem}[section]

\newtheorem{lem}[thm]{Lemma}

\newtheorem{exa}[thm]{Example}
\newtheorem{rem}[thm]{Remark}

\theoremstyle{definition}

\definecolor{wco}{rgb}{0.5,0.2,0.3}

\numberwithin{equation}{section} \theoremstyle{remark}

\def\N{\mathbb N}  
   
\def\<{\langle} \def\>{\rangle}  
    
\def\d{\text{\rm{d}}}   
  
 \def\beq{\begin{equation}}

 \def\ee{\varepsilon}

\def\L{\mathcal L}

\def\R{\mathbb{R}}

\def\L{\Lambda}
\def\W{W_0^{m,2}}
\def\H{H^m}

\def\[{{\Big[}}
\def\]{{\Big]}}

\def\({{\Big(}}
\def\){{\Big)}}

\newfam\msbfam
\font\tenmsb=msbm10 \textfont\msbfam=\tenmsb \font\sevenmsb=msbm7
\scriptfont\msbfam=\sevenmsb \font\fivemsb=msbm5
\scriptscriptfont\msbfam=\fivemsb

\def\th#1{\vspace{1mm}\noindent{\bf #1}\quad}

\def\proof{\vspace{1mm}\noindent{\it Proof}\quad}

\numberwithin{equation}{section}

\def \eref#1{\hbox{(\ref{#1})}}

\def\R{{\mathbb R}}

\def\bc{\begin{center}}
\def\ec{\end{center}}

\def \eref#1{\hbox{(\ref{#1})}}
\def\d#1#2{\frac{\displaystyle #1}{\displaystyle #2}}

\def\L{\Lambda}
\def\W{W_0^{m,2}}
\def\H{H^m}

\def\d{\text{\rm{d}}}
\def\<{\langle}
\def\>{\rangle}

\begin{document}
\title{{\Large \bf Well-posedness of Backward Stochastic Partial
Differential Equations with Lyapunov
Condition}
\footnote{Supported in part by NSFC (No. 11571147, 11671035, 11822106, 11831014),  NSF
of Jiangsu Province
 (No. BK20160004), the Qing Lan Project and
PAPD of Jiangsu Higher Education Institutions.  Financial support by the DFG through the CRC
1283``Taming uncertainty and profiting from randomness and low regularity in analysis, stochastics and their
applications'' is acknowledged.}}

\author{ {Wei Liu$^a$},  {Rongchan Zhu$^{b,c}$\footnote{Corresponding author: zhurongchan@126.com}}
\\
 \small  $a.$  School of Mathematics and Statistics, Jiangsu Normal University, Xuzhou 221116, China \\
 \small  $b.$ Department of Mathematics, Beijing Institute of Technology, Beijing 100081, China\\
 \small  $c.$ Department of Mathematics, University of Bielefeld, D-33615 Bielefeld, Germany}
\date{}
\maketitle
\begin{center}
\begin{minipage}{145mm}
{\bf Abstract.}  In this paper we show the
existence and uniqueness of strong solutions for a large class of backward SPDE  where the coefficients satisfy
a specific type Lyapunov condition instead of
the classical coercivity condition. Moreover, based on the generalized variational framework, we also use the local monotonicity condition
 to replace the standard monotonicity condition, which is applicable to various quasilinear and semilinear BSPDE models.

\vspace{3mm}

\vspace{1mm} {\bf Keywords:} BSDE; SPDE; locally monotone; Lyapunov condition
\end{minipage}
\end{center}

\section{Introduction}

The theory of backward stochastic differential equations (BSDEs) has received
 extensive investigations in the last few decades. BSDEs have been successfully applied
 in stochastic control theory, econometrics, mathematical finance,
nonlinear partial differential equations and so on,  see \cite{DE92, EM97,PP90, YZ99} and more references therein.
The study of backward stochastic partial differential equations (BSPDEs) could be traced back
to the works \cite{B83,P79}.  This subject arise in many applications of probability theory
 and stochastic processes, for instance in nonlinear filtering and
stochastic control theory for processes with incomplete information, as an adjoint equation of
the Duncan-Mortensen-Zakai filtration (see e.g. \cite{B83,HM02,HP91,T98,Z92,Z93}). In the
dynamic programming theory, some nonlinear BSPDEs as the backward stochastic
Hamilton-Jacobi-Bellman equations, are also introduced in the investigation of non-Markovian control
problems (see e.g \cite{EK09,P92}).
Recently, there are many papers studying  backward stochastic partial differential equations (see \cite{DTZ13,QT12,QTY12,SY09,T05,Z09} and the references therein). In \cite{TW16} a very general system of backward stochastic partial differential equations
is studied, and in \cite{QTY12,SY09} the authors concentrate on the study of the backward stochastic 2D Navier-Stokes equation (BSNSE).

The main aim of this work is to prove the existence and uniqueness of  solutions
for a large class of backward stochastic partial differential equations using the variational approach.
The variational framework has been used intensively for studying PDE and SPDE
where the coefficients satisfy the classical monotonicity and coercivity conditions. In the
case of deterministic equations, the theory of monotone operators started from the substantial
work of Minty \cite{M62, M63}, then it was studied systematically by Browder \cite{B63, B64} in order to obtain the existence of solutions for quasi-linear elliptic and parabolic partial differential
equations. We refer to the monograph \cite{B73} for more extensive exposition and references.
Concerning the stochastic equations, it was first investigated in the seminal works
of Pardoux \cite{P75} and Krylov and Rozovskii \cite{KR79}, where they adapted the monotonicity tricks
to prove the existence and uniqueness of solutions for a class of semilinear and quasilinear
SPDE.  Recently,
this framework has been substantially extended by the first named author and R\"{o}ckner  in \cite{LR10,LR13} for
more general class of SPDEs with coefficients satisfying the generalized coercivity and local monotonicity
conditions, hence many fundamental examples such as stochastic Burgers type equations
and stochastic 2D Navier-Stokes equations can be included into this framework now
(see \cite{L13,LR15} for more examples).

 In this paper we will show the
existence and uniqueness of strong solutions for a class of BSPDE where the coefficients satisfy
a specific type Lyapunov condition (we call it one-sided linear growth here) instead of
the classical coercivity condition. Based on \cite{LR10}, we also use the local monotonicity condition
here to replace the standard monotonicity condition.  This Lyapunov type condition (see (H3) below) is inspired by the recent work of
 \cite{L13} (see also the references therein), where  this type of condition is used to
investiage stochastic tamed 3D Navier-Stokes equations and the stochastic curve shortening
flow in the plane.
 Moreover, we should remark that our main result is
also applicable to backward stochastic 2D Navier-Stokes equations, stochastic p-Laplace equations,  stochastic fast diffusion equations, stochastic Burgers type equations and
stochastic reaction-diffusion equations. We refer to Section
3 for the details.

\section{Main Result}

First we introduce our framework in detail.
Let $(H,\langle\cdot,\cdot\rangle_H)$ be a separable Hilbert space and identified with its dual space $H^*$ by the Riesz isomorphism, and let $(V,\langle\cdot,\cdot\rangle_V)$ be a Hilbert space such that it is continuously and densely embedded into $H$. Then we have the following Gelfand triple
$$V\subset H\equiv H^*\subset V^*,$$
where $V^*$ is the dual space of $V$ (w.r.t. $\langle\cdot,\cdot\rangle_H$).

Let $(\Omega,\mathcal{F},\mathcal{F}_t,P)$ be a complete filtrated probability space, on which  a cylindrical Wiener process $\{W_t\}_{t\geq0}$ is defined on a separable Hilbert space $(U,\langle\cdot,\cdot\rangle_U)$, whose natural augmented
filtration is denoted by $\{\mathcal{F}_t, t\in[0, T]\}$ and $(L_2(U,V),\|\cdot\|_{L_2(U,V)})$ denotes the space of all Hilbert-Schmidt operators from $U$ to $V$.
We denote by $\mathcal{P}$ the $\sigma$-algebra of the
predictable sets on
$\Omega\times [0, T]$ associated with $\{\mathcal{F}_t\}_{t\geq0}$.
 For any Banach space $\mathbb{B}$, let $L_{\mathcal{F}}^p(\Omega;L^r([0,T];\mathbb{B})), p,r\in[1,\infty]$ be the set of all predictable $\mathbb{B}$-valued processes in $L^p(\Omega;L^r([0,T];\mathbb{B}))$ and let $L^p_\mathcal{F}(\Omega\times [0,T];\mathbb{B})$ to denote all predictable $\mathbb{B}$-valued processes in $L^p(\Omega\times [0,T];\mathbb{B})$.  We also use $L^p_{\mathcal{F}_T}(\Omega;\mathbb{B}), p\in[1,\infty]$ to denote the set of all $\mathcal{F}_T$-measurable random variable in $L^p(\Omega;\mathbb{B})$.

We consider the following backward stochastic partial differential equation
\begin{equation}dX_t=-A(t,X_t,Z_t)dt+Z_tdW_t, \quad t\in[0,T],\quad X(T)=\xi,\end{equation}
where $A:[0,T]\times V\times L_2(U,H)\times \Omega\rightarrow V^*$ and for any $(v,z)\in V\times L_2(U,H)$, $A(\cdot,v,z,\cdot)$ is predictable and $V^*$-valued process.

We need to suppose the following assumptions concerning the Gelfand triple.

(H0) There exists an orthogonal set $\{e_1,e_2,...\}$ in $(V,\langle\cdot,\cdot\rangle_V)$ such that it constitute an orthonormal basis of $(H,\langle\cdot,\cdot\rangle_H)$.

Suppose that there exist constants $\varepsilon\in(0,\frac{1}{2})$, $K$ and a positive adapted process $f\in L^\infty_\mathcal{F}(\Omega,L^1([0,T]))$ such that the following conditions hold for all $v,v_1,v_2\in V, \phi,\phi_1,\phi_2\in L_2(U,V)$ and $(t,\omega)\in [0,T]\times \Omega$:

(H1) (Hemicontinuity) The map $s\mapsto { }_{V^*}\!\langle A(t,v_1+sv_2,\phi),v\rangle_V$ is continuous on $\mathbb{R}$.

(H2) (Local monotonicity) There exists a locally bounded measurable function $\rho: V\rightarrow[0,+\infty)$ such that
$${ }_{V^*}\!\langle A(t,v_1,\phi_1)-A(t,v_2,\phi_2),v_1-v_2\rangle_V\leq \rho(v_2)[\|v_1-v_2\|_H^2+\|v_1-v_2\|_H\|\phi_1-\phi_2\|_{L_2(U,H)}]. $$

(H3) (One-sided linear growth) For any $n\in\mathbb{N}$, the operator $A$ maps $H_n:=\textrm{span}\{e_1,...,e_n\}$ into $V$ such that for $v\in H^n$
$$\langle A(t,v,\phi),v\rangle_V\leq f_t+\varepsilon\|\phi\|_{L_2(U,V)}^2+K\|v\|_V^2.$$

(H4) (Growth) $$\|A(t,v,\phi)\|_{V^*}\leq f_t^{1/2}+\rho(v)+K\|\phi\|_{L_2(U,V)}.$$
\vskip.10in

\th{Definition 2.1} \emph{For $\xi\in L^\infty_{\mathcal{F}_T}(\Omega,V)$ we say that $(X,Z)$ is a solution to (2.1) if
$$X\in L^2_{\mathcal{F}}(\Omega;L^\infty([0,T];V))\cap L_{\mathcal{F}}^2(\Omega;C([0,T];H)),$$
$$Z\in L^2_{\mathcal{F}}(\Omega;L^2([0,T];L_2(U,H))),$$
$$X_t=\xi+\int_t^TA(t,X_s,Z_s)ds-\int_t^TZ_sdW_s, \textrm{ in } V^*\quad P-a.s.. $$}
\vskip.10in

Now we state the main result of this work.

\th{Theorem 2.2} \emph{Suppose (H0)-(H4) hold. For any $\xi\in L^\infty_{\mathcal{F}_T}(\Omega;V)$, BSPDE (2.1) admits a unique adapted solution $(X,Z)\in L^\infty_\mathcal{F}(\Omega\times[0,T];V)\times L^2_\mathcal{F}(\Omega;L^2([0,T],L_2(U,H)))$.
Moreover, it satisfies that
$$\sup_{t\in[0,T]}\|X_t\|_V^2+\frac{1}{2}E\int_0^T\|Z_s\|_{L_2(U,V)}^2ds\leq C(\|f\|_{L_{\mathcal{F}}^\infty(\Omega;L^1([0,T]))}+\|\xi\|^2_{L_{\mathcal{F}_T}^\infty(\Omega;V)}), \quad a.s..$$}

\th{Remark 2.3} \emph{(1) In  the theorem above we assume that $\xi\in L^\infty_{\mathcal{F}_T}(\Omega,V)$, which seems quite different to the condition usually posed on the initial value of stochastic PDEs. The reason for this is that we cannot use any stopping time argument  for BSPDEs, therefore, here we have to use stochastic Gronwall-Bellman inequality to deduce an uniform estimate (see (2.5) below) to control the nonlinear term. This is also one of the main differences between BSPDEs and (standard/forward) SPDEs.}

\emph{(2) Note that (H1) and (H2) imply that $A(t,v,z)$ is locally Lipschitz continuous with respect to $z$ in the following sense:
$$\|A(t,v,z_1)-A(t,v,z_2)\|_{V^*}\leq C(\|v\|_V)\|z_1-z_2\|_{L_2(U,H)},$$
for all $t\in[0,T],\omega\in\Omega, v\in V, z_1,z_2\in L_2(U,H)$.}
\vskip.10in

The rest part of this section is devoted to the proof of main result, and we need several lemmas for this purpose.

Recall that $\{e_1,e_2,...\}\subset V$ is an orthonormal basis of $H$ and $H_n=\textrm{span}\{e_1,...,e_n\}$. Let $P_n:V^*\rightarrow H_n$ be defined by$$P_ny=\sum_{i=1}^n{ }_{V^*}\!\langle y,e_i\rangle_Ve_i,\quad y\in V^*.$$
Hence we have
$${ }_{V^*}\!\langle P_nA(t,u,z),v\rangle_V=\langle P_n A(t,u,z),v\rangle_H={ }_{V^*}\!\langle A(t,u,z),v\rangle_V,\quad u\in V,v\in H_n,z\in L_2(U,H).$$
By (H0) we have
$$\langle P_nA(t,u),v\rangle_V=\langle A(t,u),v\rangle_V.$$

 Now we consider the following projected approximation:
\begin{equation}X_t^N=P_N\xi+\int_t^TP_NA(t,X_t^N,Z_t^N)dt-\int_t^T Z_t^N dW_t,\end{equation}
where  $Z^N\in L(U,H_N)$ can be extended to a element in $L_2(U,H)$ (still denoted by $Z^N$) by setting $Z^N(e_j)=0, j\geq N+1$.

To solve  (2.2) we shall make use of the result in \cite{BDHPS03}. We fixed the filtration generated by the cylindrical Wiener process. We do not approximate $W$ by  its finite dimensional projection.  Since $W$ is a cylindrical Wiener process, we cannot apply the results in \cite[Theorem 4.2]{BDHPS03} directly.
We need to use the lemma below for the following type BSDE
\begin{equation}Y_t=\zeta+\int_t^Tg(s,Y_s,q_s)ds-\int_t^Tq_sd{W}_s,\end{equation}
where $\zeta$ is an $\mathbb{R}^N$-valued $\mathcal{F}_T$-measurable random vector, the random function $g:\Omega\times [0,T]\times \mathbb{R}^N\times L_2(U;\mathbb{R}^N)\rightarrow\mathbb{R}^N$ is $\mathcal{P}\times \mathcal{B}(\mathbb{R}^N)\times \mathcal{B}(L_2(U;\mathbb{R}^N))$-measurable.

\vskip.10in

\th{Lemma 2.4} \emph{Assume that $g$ and $\zeta$ satisfy the following four conditions:}

\emph{(C1) For some $p>1$ we have
$$E[|\zeta|^p+(\int_0^T|g(t,0,0)|dt)^p]<\infty.$$}

\emph{(C2) There exist constants $\alpha\geq0$ and $\mu\in\mathbb{R}$ such that almost surely we have for each $t\in[0,T], y,y^\prime\in\mathbb{R}^N, z,z^\prime\in L_2(U;\mathbb{R}^N)$,
$$|g(t,y,z)-g(t,y,z^\prime)|\leq \alpha\|z-z^\prime\|_{L_2(U;\mathbb{R}^N)},$$
$$\langle y-y^\prime,g(t,y,z)-g(t,y^\prime,z)\rangle\leq \mu|y-y^\prime|^2.$$}

\emph{(C3) The function $y\mapsto g(t,y,z)$ is continuous for every $(t,z)\in[0,T]\times L_2(U;\mathbb{R}^N)$.}

\emph{(C4) For any $r>0$, the stochastic process
$$\{\psi_r(t):=\sup_{|y|\leq r}|g(t,y,0)-g(t,0,0)|,t\in[0,T]\}$$
lies in the space $L^1_\mathcal{F}(\Omega\times[0,T])$.}

\emph{Then BSDE (2.3) admits a unique solution
 $$ (Y,q)\in L_\mathcal{F}^p(\Omega; C([0,T];\mathbb{R}^N))\times L^p_\mathcal{F}(\Omega;L^2([0,T];L_2(U;\mathbb{R}^N))).$$}
\vskip.10in

By the martingale representation theorem in infinite dimensional case in \cite{HP91}, we could prove Lemma 2.4. The method to prove it is standard and is a slight modification of the proof of  \cite[Theorem 4.2]{BDHPS03}, so we
omit it here. For more details we refer to \cite{Z12a, Z12b}.  The following lemma comes from \cite[Lemma 4.2]{QTY12}.
\vskip.10in
\th{Lemma 2.5} \emph{For any $M,N\in \mathbb{N}$, define
 $\varphi_n(z)=\frac{zn}{(\|z\|_{L_2(U,H_N)}\vee n)}, z\in L_2(U,H_N)$ and set
 $$A^{N,M,n}(t,y,z):=R_M(\|y\|_V)\frac{n}{h_M(t)\vee n}P_NA(t,y,\varphi_n(z)),$$
 where $R_M:\mathbb{R}\rightarrow[0,1]$ is a smooth function satisfying $R_M(r)=1, |r|\leq M$, $R_M(r)=0, |r|>M+1$, $|R_M'|\leq 1$ and
$$h_M(t):=f_t^{1/2}+\sup_{\|v\|_V\leq M}\rho(v)\in L^1(\Omega\times[0,T]). $$
Then under (H0)-(H4) $A^{N,M,n}$ satisfies the conditions (C2)-(C4) of Lemma 2.4.}

\proof Now we verify that there is a uniform constant $C_{N,M,n}>0$ such that
$$\langle A^{N,M,n}(t,X_1,Z)-A^{N,M,n}(t,X_2,Z),X_1-X_2\rangle_H\leq C_{N,M,n}\|X_1-X_2\|_H^2, \quad a.s.,$$
for $X_1,X_2\in H_N, Z\in L_2(U,H_N)$ and all $t\in[0,T]$. It holds trivially if $\|X_1\|_V>M+1$ and $\|X_2\|_V>M+1$. Thus, it is sufficient to consider the case of $\|X_2\|_V\leq M+1$. We have
$$\aligned&\langle A^{N,M,n}(t,X_1,Z)-A^{N,M,n}(t,X_2,Z),X_1-X_2\rangle
\\=&R_M(\|X_1\|_V)\frac{n}{h_M(t)\vee n}\langle P_NA(t,X_1,\varphi_n(Z))-P_NA(t,X_2,\varphi_n(Z)), X_1-X_2\rangle\\&+\frac{n}{h_M(t)\vee n}(R_M(\|X_1\|_V)-R_M(\|X_2\|_V))\langle P_NA(t,X_2,\varphi_n(Z)), X_1-X_2\rangle
\\\leq& C_{N,M,n}\|X_1-X_2\|^2, \endaligned$$
where we used (H2) and (H4) in the last inequality. The other conditions are satisfied obviously also by (H2) and (H4). $\hfill\Box$
\vskip.10in
We now recall the stochastic Gronwall-Bellman inequality from \cite[Corollary B1]{DE92}. Let $(\Omega, \mathcal{F}, \mathcal{F}_t, P)$ be a filtered probability
space whose filtration $\mathcal{F}= \{\mathcal{F}_t : t\in[0, T]\}$ satisfies the usual conditions.
Suppose that $\{Y_s\}$ and $\{X_s\}$ are optional integrable processes and $\alpha$ is a nonnegative
constant. If for all $t$, the map $s\rightarrow E[Y_s|\mathcal{F}_t]$ is continuous almost surely and
$$Y_t\leq E[\int_t^T(X_s + \alpha Y_s) ds + Y_T|\mathcal{F}_t],$$
then we have almost surely
$$Y_t\leq e^{\alpha(T-t)}E[Y_T |\mathcal{F}_t] + E[\int_t^Te^{\alpha(s-t)}X_sds|\mathcal{F}_t], t\in[0,T].$$

\vskip.10in
\th{Lemma 2.6} \emph{Suppose that Assumptions (H1)-(H4) hold. For any $\xi\in L^\infty_{\mathcal{F}_T}(\Omega;V)$, the projected problem (2.2) admits a unique adapted solution $$(X^N,Z^N)\in L^2_\mathcal{F}(\Omega;C([0,T];H))\times L_\mathcal{F}^2(\Omega;L^2([0,T],L_2(U,H))).$$ }

\proof [Existence]. By Lemmas 2.4 and 2.5 there exists a unique solution $(X^{N,M,n}, Z^{N,M,n})$ to the following BSDE
$$X^{N,M,n}_t=\xi^N+\int_t^TA^{N,M,n}(s,X^{N,M,n}_s,Z^{N,M,n}_s)ds-\int_t^TZ^{N,M,n}_sdW_s,$$
for $\xi^N=P_N\xi$ and
$$X^{N,M,n}\in L^2_{\mathcal{F}}(\Omega;L^\infty([0,T];V))\cap L_{\mathcal{F}}^2(\Omega;C([0,T];H)),$$
$$Z^{N,M,n}\in L^2_{\mathcal{F}}(\Omega;L^2([0,T];L_2(U,H))).$$
Now by It\^{o}'s formula and using (H3) we have
\begin{equation}\aligned\|X^{N,M,n}_t\|_V^2=&\|\xi^N\|_V^2+2\int_t^T\langle A^{N,M,n}(s,X^{N,M,n}_s,Z^{N,M,n}_s),X_s^{N,M,n}\rangle_Vds
\\&-\int_t^T\|Z_s^{N,M,n}\|_{L_2(U,V)}^2ds -2\int_t^T\langle X_s^{N,M,n},Z_s^{N,M,n}dW_s\rangle_V
\\\leq&\|\xi^N\|_V^2+2\int_t^T(f_s+\varepsilon\|Z_s^{N,M,n}\|_{L_2(U,V)}^2+K\|X^{N,M,n}_s\|_V^2)ds
\\&-\int_t^T\|Z_s^{N,M,n}\|_{L_2(U,V)}^2ds-2\int_t^T\langle X_s^{N,M,n},Z_s^{N,M,n}dW_s\rangle_V.\endaligned\end{equation}
By the Burkholder-Davis-Gundy inequality we have
$$\aligned  & E\left[\sup_{\tau\in[t,T]}\left|\int_\tau^T\langle X_s^{N,M,n},Z_s^{N,M,n}dW_s\rangle_V\right|\right]\\
\leq &2E\left[\sup_{\tau\in[t,T]}\left|\int_t^\tau\langle X_s^{N,M,n},Z_s^{N,M,n}dW_s\rangle_V\right|\right]\\\leq &CE\left[\int_t^T\| X_s^{N,M,n}\|_V^2\|Z_s^{N,M,n}\|_{L_2(U,V)}^2ds\right]^{1/2}
\\\leq&C(N)\left\{E\sup_{s\in[t,T]}\| X_s^{N,M,n}\|_V^2+E\int_t^T\|Z_s^{N,M,n}\|_{L_2(U,H)}^2ds \right\},\endaligned$$
where we used $\|Z_s^{N,M,n}\|_{L_2(U,V)}\leq C(N)\|Z_s^{N,M,n}\|_{L_2(U,H)} $ in the last step.
Taking conditional expectation on both sides of (2.4) we obtain
$$\aligned\|X^{N,M,n}_t\|_V^2+\frac{1}{2}E_{\mathcal{F}_t}\int_t^T\|Z_s^{N,M,n}\|_{L_2(U,V)}^2ds\leq&E_{\mathcal{F}_t}\|\xi^N\|_V^2+2E_{\mathcal{F}_t}\int_t^T(f_s+K\|X^{N,M,n}_s\|_V^2)ds.\endaligned$$
By the stochastic Gronwall-Bellman inequality, we conclude that
\begin{equation}\sup_{t\in[0,T]}\|X^{N,M,n}_t\|_V^2+\frac{1}{2}E\int_0^T\|Z_s^{N,M,n}\|_{L_2(U,V)}^2ds\leq C(\|f\|_{L_{\mathcal{F}}^\infty(\Omega;L^1([0,T]))}+\|\xi\|^2_{L_{\mathcal{F}_T}^\infty(\Omega;V)}), \quad a.s.,\end{equation}
where $C$ is a constant independent of $N,M,n$. Now we deduce that there exists a positive constant $K_1$ independent of $N,M$ and $n$ such that
$$\sup_{t\in[0,T]}\|X^{N,M,n}_t\|_V^2+\frac{1}{2}E\int_0^T\|Z_s^{N,M,n}\|_{L_2(U,V)}^2ds\leq K_1, \quad a.s..$$
Then letting $M=K_1+1$ be fixed, we have $R_M(\|X^{N,M,n}_s\|_V)\equiv1$. Now we write  $(X^{N,n},Z^{N,n})$ instead of $(X^{N,M,n},Z^{N,M,n})$ below.
Then  there exists a positive constant $K_2$ independent of $N,M$ and $n$ such that$$\|A(t,X^{N,n}_t,z_1)-A(t,X^{N,n}_t,z_2)\|_{V^*}\leq K_2\|z_1-z_2\|_{L_2(U,H)},$$
$$\rho(X^{N,n}_t)+\rho^2(X^{N,n}_t)\leq K_2.$$
For $j\in\mathbb{N}$, set $({\bar{X}}^N,\bar{Z}^N)=(X^{N,n+j}-X^{N,n},Z^{N,n+j}-Z^{N,n})$. Applying  It\^{o}'s formula we get
$$\aligned &e^{K_2t}\|\bar{X}^N_t\|_H^2+\int_t^Te^{K_2s}\|\bar{Z}_s^{N}\|_{L_2(U,H)}^2ds\\\leq& 2\int_t^Te^{K_2s}\langle A^{N,n+j}(s,X^{N,n+j}_s,Z^{N,n+j}_s)-A^{N,n}(s,X^{N,n}_s,Z^{N,n}_s),\bar{X}_s^{N}\rangle_H ds\\&-2\int_t^Te^{K_2s}\langle \bar{X}_s^{N},\bar{Z}_s^{N}dW_s\rangle_H-K_2\int_t^Te^{K_2s}\|\bar{X}^N_s\|_H^2ds\\\leq& 2\int_t^Te^{K_2s}\langle A^{N,n+j}(s,X^{N,n+j}_s,Z^{N,n+j}_s)-A^{N,n+j}(s,X^{N,n}_s,Z^{N,n}_s),\bar{X}_s^{N}\rangle_H ds\\&+2\int_t^Te^{K_2s}\langle A^{N,n+j}(s,X^{N,n}_s,Z^{N,n}_s)-A^{N,n}(s,X^{N,n}_s,Z^{N,n}_s),\bar{X}_s^{N}\rangle_H ds\\&-2\int_t^Te^{K_2s}\langle \bar{X}_s^{N},\bar{Z}_s^{N}dW_s\rangle_H-K_2\int_t^Te^{K_2s}\|\bar{X}^N_s\|_H^2ds\\\leq& 2\int_t^Te^{K_2s}\langle \frac{n+j}{h_{K_1+1}(t)\vee (n+j)} P_NA(t,X^{N,n},\varphi_{n+j}(Z^{N,n}))\\&-\frac{n}{h_{K_1+1}(t)\vee n}P_NA(t,X^{N,n},\varphi_n(Z^{N,n})), \bar{X}_s^{N}\rangle_H ds-2\int_t^Te^{K_2s}\langle \bar{X}_s^{N},\bar{Z}_s^{N}dW_s\rangle_H\\\leq& 8CK_1K_2\int_t^Te^{K_2s}\left[\|Z^{N,n}_s\|_{L_2(U,H)}1_{\{\|Z^{N,n}_s\|_{L_2(U,H)}>n\}}
+\|Z^{N,n}_s\|_{L_2(U,V)}1_{\{h_{K_1+1}>n\}}+2h_{K_1+1}1_{\{h_{K_1+1}>n\}}\right] ds\\&-2\int_t^Te^{K_2s}\langle \bar{X}_s^{N},\bar{Z}_s^{N}dW_s\rangle_H,\endaligned$$
where we used $\|Z^{N,n}\|_{L_2(U,H)}\sim\|Z^{N,n}\|_{L_2(U,V)}$.

On the other hand by the BDG inequality we have
$$E[\sup_{\tau\in[t,T]}|\int_\tau^Te^{K_2s}\langle \bar{X}_s^{N},\bar{Z}_s^{N}dW_s\rangle|]\leq \varepsilon_0E[\sup_{s\in[t,T]}(e^{K_2s}\|\bar{X}_s^{N}\|^2)]+CE[\int_t^Te^{K_2s}\|\bar{Z}_s^{N}\|_{L_2(U,H)}^2ds]. $$
By (2.5) and $h_{K_1+1}\in L^1(\Omega\times[0,T])$ we conclude that $(X^{N,n},Z^{N,n})$ is a Cauchy sequence in
$L^2_\mathcal{F}(\Omega;C([0,T];H))\times L_\mathcal{F}^2(\Omega;L^2([0,T],L_2(U,H)))$. Denote the limit by $(X^N,Z^N)$. It is easily checked that $(X^N,Z^N)$ is a solution to (2.2).

\vskip.10in
[Uniqueness]. Suppose $(X_1^N,Z_1^N)$ and $(X^N_2,Z_2^N)$ are two solutions of the projected equation (2.2). In the proof of uniqueness we use $X(t)$ to denote $X_t$.
Denote $(\tilde{X}^N,{\tilde{Z}}^N):=(X_1^N-X_2^N,Z_1^N-Z_2^N)$. By the same arguments as above we obtain (2.5) also holds for $(\tilde{X}^N,{\tilde{Z}}^N)$. Define
$$r(t):=2\int_0^t\rho(X_2^N(s))+\rho^2(X_2^N(s))ds.$$
An application of It\^{o}'s formula and (H2) yields that
$$\aligned e^{r(t)}\|\tilde{X}^N(t)\|_H^2=&\int_t^Te^{r(s)}[2\langle P_NA(s,X_1^N(s),Z_1^N(s))-P_NA(s,X_2^N(s),Z_2^N(s)),\tilde{X}^N(s)\rangle_H
\\&-\|{\tilde{Z}}^N\|^2_{L_2(U,H)}-2\|\tilde{X}^N(s)\|_H^2(\rho(X_2^N(s))+\rho^2(X_2^N(s)))]ds\\&-2\int_t^Te^{r(s)}\langle\tilde{X}^N(s),{\tilde{Z}}^NdW_s \rangle_H
\\\leq&\int_t^Te^{r(s)}[2\|\tilde{X}^N(s)\|_H^2\rho(X_2^N(s))+2\|\tilde{X}^N(s)\|\rho(X_2^N(s))\|{\tilde{Z}}^N(s)\|_{L_2(U,H)}
\\&-\|{\tilde{Z}}^N\|^2_{L_2(U,H)}-2\|\tilde{X}^N(s)\|_H^2(\rho(X_2^N(s))+\rho^2(X_2^N(s)))]ds
\\&-2\int_t^Te^{r(s)}\langle\tilde{X}^N(s),{\tilde{Z}}^NdW_s\rangle_H \endaligned$$
Taking conditional expectation on both sides we have for any $t\in[0,T]$,
$$e^{r(t)}\|\tilde{X}^N(t)\|_H^2+\frac{1}{2}E_{\mathcal{F}_t}[\int_t^Te^{r(s)}\|{\tilde{Z}}^N\|^2_{L_2(U,H)}ds]\leq 0, a.s.,$$
which implies the uniqueness. $\hfill\Box$
\vskip.10in

Now we can finish the proof of Theorem 2.2.

\proof [Existence] By the same arguments as in the proof of (2.5) we obtain
\begin{equation}\sup_{t\in[0,T]}\|X^{N}_t\|_V^2+\frac{1}{2}E\int_0^T\|Z_s^{N}\|_{L_2(U,V)}^2ds\leq C(\|f\|_{L_{\mathcal{F}}^\infty(\Omega;L^1([0,T]))}+\|\xi\|^2_{L_{\mathcal{F}_T}^\infty(\Omega;V)}), \quad a.s.,\end{equation}
where $C$ is independent of $N$.
By (H4) we have
$$E\int_0^T\|A(t,X^N,Z^N)\|_{V^*}^2dt\leq C.$$
Then there exists a subsequence $N_k\rightarrow\infty$ such that

(i) $X^{N_k}\rightarrow\bar{X}$ weakly in $L^2_\mathcal{F}(\Omega;L^2([0,T],V))$ and weakly star in $L^\infty_\mathcal{F}((\Omega\times [0,T]),V)$.

(ii) $Y^{N_k}:=A(t,X^{N_k},Z^{N_k})\rightarrow Y$  weakly in $L^2_\mathcal{F}(\Omega;L^2([0,T],V^*))$.

(iii) $Z^{N_k}\rightarrow Z$ weakly in $L^2_\mathcal{F}(\Omega;L^2([0,T],L_2(U,V)))$ and hence
$$\int_\cdot^TZ^{N_k}_sdW_s\rightarrow\int_\cdot^TZ_sdW_s$$
weakly in $L^\infty([0,T];L^2(\Omega,H))$.
Now we define the following process
$$X_t:=\xi+\int_t^TY_sds-\int_t^TZ_sdW_s, t\in[0,T], $$
then it is easy to show that $X=\bar{X},  dt\times P$-a.e.
By \cite[Theorem 4.2.5]{LR15},  we conclude that $X\in L^2_\mathcal{F}(\Omega;C([0,T],H))$ and by (2.6) we obtain
$$\sup_{t\in[0,T]}\|X_t\|_V^2+\frac{1}{2}E\int_0^T\|Z_s\|_{L_2(U,V)}^2ds\leq C(\|f\|_{L_{\mathcal{F}}^\infty(\Omega;L^1([0,T]))}+\|\xi\|^2_{L_{\mathcal{F}_T}^\infty(\Omega;V)}), \quad a.s.,$$

Now it is sufficient to show that
$$A(\cdot,\bar{X},Z)=Y, dt\times P-a.e.. $$

For $v\in L_\mathcal{F}^\infty(\Omega\times [0,T];V)$ we define $$r_t:=\int_0^t\rho(v_s)+\rho^2(v_s)ds.$$
Applying the It\^{o}'s formula we have
$$\aligned &E[e^{r_t}\|X^N_t\|_H^2-e^{r_T}\|X^N_T\|_H^2]\\=&E[\int_t^Te^{r_s}(2\langle P_NA(s,X^N_s,Z^N_s), X^N_s\rangle_H-\|Z^N_s\|_{L_2(U,H)}^2\\&-(\rho(v_s)+\rho^2(v_s))\|X^N_s\|_H^2)ds]\\=&E[\int_t^Te^{r_s}(2{ }_{V^*}\!\langle A(s,X^N_s,Z^N_s), X^N_s\rangle_V-\|Z^N_s\|_{L_2(U,H)}^2\\&-(\rho(v_s)+\rho^2(v_s))\|X^N_s\|_H^2)ds]\\=&E[\int_t^Te^{r_s}(2{ }_{V^*}\!\langle A(s,X^N_s,Z^N_s)-A(s,v_s,Z_s), X^N_s-v_s\rangle_V\\&-\|Z^N_s-Z_s\|_{L_2(U,H)}^2-(\rho(v_s)+\rho^2(v_s))\|X^N_s-v_s\|_H^2)ds]\\&+E[\int_t^Te^{r_s}(2{ }_{V^*}\!\langle A(s,X^N_s,Z^N_s)-A(s,v_s,Z_s), v_s\rangle_V\\&+2{ }_{V^*}\!\langle A(s,v_s,Z_s), X^N_s\rangle_V-2\langle Z^N_s,Z_s\rangle_{L_2(U,H)}+\|Z_s\|_{L_2(U,H)}\\&-(\rho(v_s)+\rho^2(v_s))2\langle X^N_s,v_s\rangle_H-\|v_s\|_H^2)ds]\\\leq&E[\int_t^Te^{r_s}(2{ }_{V^*}\!\langle A(s,X^N_s,Z^N_s)-A(s,v_s,Z_s), v_s\rangle_V\\&+2{ }_{V^*}\!\langle A(s,v_s,Z_s), X^N_s\rangle_V-2\langle Z^N_s,Z_s\rangle_{L_2(U,H)}+\|Z_s\|_{L_2(U,H)}^2\\&-(\rho(v_s)+\rho^2(v_s))(2\langle X^N_s,v_s\rangle_H-\|v_s\|_H^2))ds].\endaligned$$
Letting $N\rightarrow\infty$, by (H2) and the lower semicontinuity, we have for any nonnegative $\psi\in L^\infty([0,T])$,
\begin{equation}\aligned &E[\int_0^T\psi_t(e^{r_t}\|X_t\|_H^2-e^{r_T}\|X_T\|_H^2)dt]\\\leq&E\int_0^T[\int_t^Te^{r_s}(2{ }_{V^*}\!\langle Y_s-A(s,v_s,Z_s), v_s\rangle_V\\&+2{ }_{V^*}\!\langle A(s,v_s,Z_s), X_s\rangle_V-2\langle Z_s,Z_s\rangle_{L_2(U,H)}+\|Z_s\|_{L_2(U,H)}^2\\&-(\rho(v_s)+\rho^2(v_s))(2\langle X_s,v_s\rangle_H-\|v_s\|_H^2))ds]dt. \endaligned\end{equation}
By  It\^{o}'s formula we have
\begin{equation}\aligned &E[e^{r_t}\|X_t\|_H^2-e^{r_T}\|X_T\|_H^2]\\=&E[\int_t^Te^{r_s}(2{ }_{V^*}\!\langle Y_s, X_s\rangle_V-\|Z_s\|_{L_2(U,H)}^2-(\rho(v_s)+\rho^2(v_s))\|X_s\|_H^2)ds].\endaligned\end{equation}
Combining (2.7) with (2.8) we obtain that
\begin{equation}\aligned &E[\int_0^T\psi_t\int_t^Te^{r_s}(2{ }_{V^*}\!\langle Y_s-A(s,v_s,Z_s), X_s-v_s\rangle_V-(\rho(v_s)+\rho^2(v_s))(\|X_s-v_s\|_H^2))dsdt]\leq0. \endaligned\end{equation}
Taking $v=X-\varepsilon\phi w$ for $\varepsilon>0$ and $\phi\in L_\mathcal{F}^\infty(\Omega\times [0,T];dt\times P;\mathbb{R})$ and $w\in V$. Then we divide by $\varepsilon$ and letting $\varepsilon\rightarrow0$ to derive that
\begin{equation}\aligned &E[\int_0^T\psi_t\int_t^Te^{r_s}2{ }_{V^*}\!\langle Y_s-A(s,X_s,Z_s), w\rangle_V]\leq0. \endaligned\end{equation}
Then $Y=A(\cdot,X,Z)$ follows from the arbitrariness of $\psi$ and $w$.

[Uniqueness] Suppose that $(X_1,Z_1)$ and $(X_2,Z_2)$ are two solutions of the  problem (2.1). In the proof of uniqueness we use $X(t)$ to denote $X_t$.
Denote $(\tilde{X},{\tilde{Z}}):=(X_1-X_2,Z_1-Z_2)$. Define
$$r_1(t):=2\int_0^t\rho(X_2(s))+\rho^2(X_2(s))ds.$$
An application of It\^{o}'s formula yields that
$$\aligned e^{r_1(t)}\|\tilde{X}(t)\|_H^2=&\int_t^Te^{r_1(s)}[2{ }_{V^*}\!\langle A(s,X_1(s),Z_1(s))-A(s,X_2(s),Z_2(s)),\tilde{X}(s)\rangle_V
\\&-\|{\tilde{Z}}\|^2_{L_2(U,H)}-2\|\tilde{X}(s)\|_H^2(\rho(X_2(s))+\rho^2(X_2(s)))]ds\\&-2\int_t^Te^{r_1(s)}\langle\tilde{X}(s),{\tilde{Z}}dW_s \rangle_H
\\\leq&\int_t^Te^{r_1(s)}[2\|\tilde{X}(s)\|_H^2\rho(X_2(s))+2\|\tilde{X}(s)\|\rho(X_2(s))\|{\tilde{Z}}(s)\|_{L_2(U,H)}
\\&-\|{\tilde{Z}}\|^2_{L_2(U,H)}-2\|\tilde{X}(s)\|_H^2(\rho(X_2(s))+\rho^2(X_2(s)))]ds
\\&-2\int_t^Te^{r(s)}\langle\tilde{X}(s),{\tilde{Z}}dW_s\rangle_H \endaligned$$
Taking conditional expectations on both sides we have for any $t\in[0,T]$,
$$e^{r_1(t)}\|\tilde{X}(t)\|_H^2+\frac{1}{2}E_{\mathcal{F}_t}[\int_t^Te^{r_1(s)}\|{\tilde{Z}}\|^2_{L_2(U,H)}ds]\leq 0, a.s.,$$
which implies the uniqueness.
 $\hfill\Box$

\section{Applications}

 Let $\Lambda$ be a  bounded domain in $\mathbb{R}^d$ with sufficiently smooth boundary and
$C_0^\infty(\Lambda, \R^d)$ denote the set of all smooth functions from $\L$ to $\R^d$
with compact support.  For $p\ge 1$, let $\left(L^p(\L, \R^d), \|\cdot\|_{L^p} \right) $ be the vector valued $L^p$-space.
 For any integer $m\ge 0$, let $W_0^{m,2}$ denote the standard  Sobolev space on $\L$
with values in $\R^d$,  $i.e.$  the closure of $C_0^\infty(\Lambda, \R^d)$ with respect to the following  norm:
$$ \|u\|_{\W}^2 = \left( \sum_{0\le |\alpha|\le m} \int_{\L} |D^\alpha u|^2  \d x \right)^2.  $$

For the reader's convenience, we recall the following Gagliardo-Nirenberg interpolation inequality,
 which is used very often
in the study of PDE theory.

\begin{lem}
If $q\in [1,\infty]$ such that
$$    \frac{1}{q}= \frac{1}{2} - \frac{m \gamma}{d},\ 0\le  \gamma\le 1,  $$
then there exists a constant $C_{m,q}>0$ such that  for any $u\in \W$,
\begin{equation}\label{GN inequality}
 \|u\|_{L^q} \le C_{m,q} \|u\|_{\W}^\gamma \|u\|_{L^2}^{1-\gamma}.
\end{equation}
\end{lem}

Now we define
$$  H^{m}:=\left\{ u\in\W: \  div(u)=0     \right\}.$$
The norm of $\W$ restricted to $H^m$ will be denoted by $\|\cdot\|_{\H}$. Note that
$H^0$ is a closed linear subspace of the Hilbert space $L^2(\L, \R^3)$.

For all  the examples in  below,  $\{W_t\}_{t\geq0}$ denotes a cylindrical Wiener process on a
separable Hilbert space $(U, \<\cdot, \cdot\>_U )$  $w.r.t$  a complete filtered probability space
$(\Omega,\mathcal{F},\mathcal{F}_t,\mathbb{P})$.

\subsection{Backward stochastic (generalized) curve shortening flow and backward singular stochastic $p$-Laplace equations}
The study of the motion by mean curvature of curves and surfaces attracts more and more attentions in recent years.  It
not only connects to many  interesting mathematical theories such as nonlinear PDEs, geometric measure theory, asymptotic analysis and
singular perturbations, but also has important applications in image processing and  materials science etc (cf.\cite{S97,Zhu02}).
The incorporation of stochastic perturbations has also been  widely used in these models, where the noise can come from the
thermal fluctuations, impurities and the atomistic processes describing the surface motions. However, the mathematical theory for
the study of those stochastic models are quite incomplete (cf.\cite{ER11} and the references therein).

The stochastic curve shortening flow (cf.\cite{ER11,ERS11}) is  formulated in the  following form:
\begin{equation*}
dX_t =  \frac{\partial_x^2 X_t}{1+(\partial_x X_t)^2 }  dt + \sigma(X_t) d W_t,
\end{equation*}
where $\partial_x, \partial_x^2$ denote the first and second (spatial) derivative, and $\sigma$ satisfies some suitable conditions.

The deterministic part is a simplified model in geometric PDE theory which describes the motion by mean curvature of embedded surfaces
 (in the present model
the surface is just  some curve  in the 2-dimensional plane), we refer to \cite{ER11} for more detailed exposition on the model.
The random forcing was introduced to refine the model by taking the influence of thermal noise into account.

Based on the crucial observation
$$ \frac{\partial_x^2 X_t}{1+(\partial_x  X_t)^2 } = \partial_x \left(  \arctan (\partial_x X_t) \right),  $$
this equation has been investigated in \cite{ER11,ERS11}  using the variational framework with following Gelfand triple:
$$    V:= W_0^{1,2} ([0,1]) \subseteq H:=  L^{2} ([0,1]) \subseteq V^*= W^{-1,2} ([0,1]).     $$

The first example here  is the equation of backward stochastic curve shortening flow, and we consider the following form of BSPDE, which covers a large class of stochastic evolution equations such as stochastic curve shortening flow (with some nonlinear perturbations), stochastic $p$-Laplace equations and stochastic reaction-diffusion equations. For simplicity we only formulate the result
for 1 dimensional underlying domain $[0, 1]$ here.
\begin{equation}\label{CSF}
dX_t = -\left[ \partial_x \left( \bar{f} (\partial_x X_t) \right) +g(X_t) +h (t,X_t,Z_t)  \right] dt + Z_t d W_t,\  X_T=\xi.
\end{equation}

\begin{exa}
Suppose  that  functions  $\bar{f}, g\in C^1(\R)$ and  there exist  constants $C, p\ge 2$ such that
\begin{equation}\begin{split}\label{F}
&  \bar{f}^\prime(x) \ge 0,  \  |\bar{f}(x)| \le C(1+|x|), \  x\in\R; \\
& g^\prime (x)\le C, \  |g(x)| \le C(1 +|x|^{p-1}), \    x\in \R;\\
& \left(g(x)-g(y)\right) \left(x-y \right) \le C(1+|y|^p) |x-y|^2, \ x,y\in\R,
\end{split}
\end{equation}
and $h:[0,T]\times V\times L_2(U,H)\times \Omega\rightarrow V^*$ satisfies (H0)-(H4).
Then for any $\xi\in L^\infty_{\mathcal{F}_T}(\Omega;V)$,
\eref{CSF} admits a unique adapted solution $(X,Z)\in L^\infty_\mathcal{F}(\Omega\times[0,T];V)\times L^2_\mathcal{F}(\Omega;L^2([0,T],L_2(U,H)))$.
Moreover, it satisfies that
$$\sup_{t\in[0,T]}\|X_t\|_V^2+\frac{1}{2}E\int_0^T\|Z_s\|_{L_2(U,V)}^2ds\leq C(1+\|f\|_{L_{\mathcal{F}}^\infty(\Omega;L^1([0,T]))}+\|\xi\|^2_{L_{\mathcal{F}_T}^\infty(\Omega;V)}), \quad a.s..$$
\end{exa}

\begin{proof} We consider the
following Gelfand triple:
$$    V:= W_0^{1,2} ([0,1]) \subseteq H:=  L^{2} ([0,1]) \subseteq V^*= W^{-1,2} ([0,1]).     $$
$(H0)$  holds since  all eigenvectors
$\{ e_i, i=1,2,\cdots \}$ of the Laplace operator  constitute  an orthonormal basis of $H$ and an orthogonal set in $V$.

By the assumptions on $\bar{f}$   we have
\begin{equation*}
 \begin{split}
&  \<\partial_x(\bar{f}(\partial_x v)), v\>_V =-\int_0^1 \bar{f}^\prime(\partial_x v)
\left( \partial_x^2 v \right)^2 dx\le 0, \  v\in H_n\subseteq V; \\
& \|\partial_x(\bar{f}(\partial_x  v))\|_{V^*} \le \| (\bar{f}(\partial_x v))\|_{H} \le C(1+\|v\|_V),\
v\in V; \\
 &  _{V^*}\<\partial_x(\bar{f}(\partial_x u))-\partial_x(\bar{f}(\partial_x u)) , u-v\>_V
=-\int_0^1 \left( \bar{f}(\partial_x u)- \bar{f}(\partial_x u) \right) \left( \partial_x u-\partial_x v \right) dx \le 0,  \ u, v\in V.
 \end{split}
\end{equation*}

We now show that $(H1)$-$(H4)$  hold for the term $g$ in the drift. By the continuity of $g$ and dominated convergence theorem it is easy to  show  that $(H1)$ holds.

By (\ref{F}) and Sobolev's  inequality we have
\begin{equation*}
 \begin{split}
 &  _{V^*}\<g(u)-g(v), u- v\>_V \\
 =& \int_0^1 (g(u)-g(v)) ( u- v) dx\\
\le &  C(1+\|v\|_{L^\infty}^p )    \int_0^1 |u-v|^2 dx \\
\le & C(1+\|v\|_V^p ) \|u-v\|_H^2,  \ u, v\in V,
 \end{split}
\end{equation*}
i.e. $(H2)$ holds with $\rho(v)=\|v\|_V^p$.

$(H3)$ also holds since  (\ref{F}) implies that
$$   \<g(v), v\>_V=-\<g(v), \partial_x^2 v\>_H=\int_0^1 g^\prime(v) (\partial_x v)^2 dx \le C \|v\|_V^2, \  v\in H_n\subseteq V.   $$

  $(H4)$  follows from the following estimate:
$$\|g(v)\|_{V^*}\le C\|g(v)\|_{L^1} \le C(1+\|v\|_{L^\infty}^{p-1} )
\le C(1+\|v\|_{V}^{p-1}) , \   v\in V.  $$

Then by the assumptions of $h$, it is easy to show that $(H1)$-$(H4)$  hold for the term $\partial_x(\bar{f}(\partial_x X_t))+g(X_t)+h(t,X_t,Z_t)$. Therefore, the conclusion follows from Theorem 2.2.
\end{proof}

\begin{rem}
 (1) If we take $\bar{f}(x)=\arctan x$ and  $g\equiv 0$, then (\ref{CSF}) reduces back to the model of
 backward stochastic curve shortening flow.

(2) The simple example of $g$ satisfying (\ref{F}) is  any polynomial of  odd degree with negative leading coefficients. Hence (\ref{CSF}) also covers backward stochastic reaction-diffusion equations (i.e. $\bar{f}(x)=x$).

(3) If $\bar{f}(x)=|x|^{p-2}x (1<p\leq 2)$, then (\ref{CSF}) covers the singular backward stochastic $p$-Laplace equations.

(3) If $\bar{f}(x)=|x|^{p-2}x (p>2)$, then (\ref{CSF}) reduces to the degenerate backward stochastic $p$-Laplace equations and the result above can not be applied to this case.
\end{rem}

\subsection{Backward stochastic fast diffusion equations}
Let $\Lambda$ be a bounded open domain in $\R^d$ with smooth boundary and $\Delta$ be the standard Laplace operator with Dirichlet boundary condition.
Stochastic fast  diffusion equations with general multiplicative noise has been studied a lot in recent years (see e.g. \cite{LR15,RRW07,KR79}). In this work,
we consider the following backward stochastic fast  diffusion equations:
 \begin{equation}\label{FDE}
 dX_t = -\left( \Delta  \Psi(X_t)+ h(t,X_t,Z_t)  \right) dt +Z_t dW_t, \ X_T=\xi,
\end{equation}
where $\Psi:  \R \to \R$  is measurable.  In
particular, if $\Psi(s)=s^r:= |s|^{r-1}s$ for
some $r\in (0,1)$, then (\ref{FDE}) reduces back to the classical
backward stochastic fast  diffusion equations.

Using the  Gelfand triple
$$    V:= L^2 (\Lambda) \subseteq H:=  W^{-1,2} (\Lambda) \subseteq V^*= (L^2 (\Lambda))^* ,    $$
we obtain the following well-posedness result for equation (\ref{FDE}).

\begin{exa}
Suppose  that  $h:[0,T]\times V\times L_2(U,H)\times \Omega\rightarrow V^*$ satisfies (H0)-(H4),  $\Psi \in C^1(\R)$ and  there exists a   constant $C>0$ such that
\begin{equation*}\begin{split}
&  \Psi^\prime(x) \ge 0,  \  |\Psi(x)| \le C(1+|x|), \  x\in\R.
\end{split}
\end{equation*}
Then for any $\xi\in L^\infty_{\mathcal{F}_T}(\Omega;V)$,
\eref{FDE} has a unique  adapted solution $(X,Z)\in L^\infty_\mathcal{F}(\Omega\times[0,T];V)\times L^2_\mathcal{F}(\Omega;L^2([0,T],L_2(U,H)))$.
In particular, we have
$$\sup_{t\in[0,T]}\|X_t\|_V^2+\frac{1}{2}E\int_0^T\|Z_s\|_{L_2(U,V)}^2ds\leq C(1+\|f\|_{L_{\mathcal{F}}^\infty(\Omega;L^1([0,T]))}+\|\xi\|^2_{L_{\mathcal{F}_T}^\infty(\Omega;V)}), \quad a.s..$$
\end{exa}
\begin{proof}
 According to the classical result for \eref{FDE} (cf. \cite[Example 4.1.11]{PR07}), here we only need to verify
the one-sided linear growth condition
$(H3)$ for (\ref{FDE}). In fact, we have
$$\<\Delta \Psi(v) , v\>_V=-\int_\Lambda \Psi^\prime(v) |\nabla v|^2 dx
 \le 0,  \ v\in H^n.   $$
Therefore, the assertions follow  directly from Theorem 2.6.
\end{proof}

\begin{rem}
(1) If $\Psi(x)=|x|^{r-1}x  (0<r< 1)$,  then (\ref{FDE}) reduces to the backward stochastic fast  diffusion equations.

(2) If $\Psi(x)=|x|^{r-1}x (r>1)$, then (\ref{FDE}) is the  backward stochastic porous medium equations, and the result above can not be applied to this case.
\end{rem}

\subsection{Backward stochastic Burgers type and  reaction-diffusion equations}
The main result in this paper is also applicable to semilinear type BSPDE  which is formulated as follows:
\begin{equation}\label{RDE}
dX_t = -\left( \partial_x^2 X_t + \bar{f}(X_t)\partial_x X_t+  g (X_t)+h(t,X_t,Z_t)\right) dt + Z_t d W_t,\  \ X_T=\xi.
\end{equation}

Consider the Gelfand triple
$$    V:= W_0^{1,2} ([0,1]) \subseteq H:=  L^{2} ([0,1]) \subseteq V^*= W^{-1,2} ([0,1]),     $$
we have the following result concerning the well-posedness of equation \eref{RDE}.

\begin{exa}
Suppose  that  $\bar{f}$  is a bounded Lipschitz function
on $\R$ and $g\in C^1(\R)$ and there exists constants $C, p\ge 2$ such that
\begin{equation*}\begin{split}\label{F1}
& \left(g(x)-g(y)\right) \left(x-y \right) \le C(1+|y|^p) |x-y|^2, \ x,y\in\R;\\
& |g(x)| \le C(1 +|x|^{p-1}), \  \  x\in \R; \\
 &   g^\prime (x)\le C, \ x\in \R,
\end{split}
\end{equation*}
and $h:[0,T]\times V\times L_2(U,H)\times \Omega\rightarrow V^*$ satisfies (H0)-(H4).
Then for any $\xi\in L^\infty_{\mathcal{F}_T}(\Omega;V)$,
\eref{RDE} has a unique  adapted solution $(X,Z)\in L^\infty_\mathcal{F}(\Omega\times[0,T];V)\times L^2_\mathcal{F}(\Omega;L^2([0,T],L_2(U,H)))$.
In particular, we have
$$\sup_{t\in[0,T]}\|X_t\|_V^2+\frac{1}{2}E\int_0^T\|Z_s\|_{L_2(U,V)}^2ds\leq C(1+\|f\|_{L_{\mathcal{F}}^\infty(\Omega;L^1([0,T]))}+\|\xi\|^2_{L_{\mathcal{F}_T}^\infty(\Omega;V)}), \quad a.s..$$
\end{exa}

\begin{proof} Combining with the result in the previous example, here we only need to show $(H1)$-$(H4)$ hold
for the term $\partial_x^2  + \bar{f}(\cdot)\partial_x$.

According to the result showed in \cite[Example 3.2]{LR10}, $(H1),(H2)$ and $(H4)$ hold.

Since $\bar{f}$ is bounded, by H\"{o}lder's inequality and Young's inequality we have
\begin{equation*}
 \begin{split}
 &  \<\partial_x^2 v+\bar{ f}(v)\partial_x v,  v\>_V \\
=& - \<\partial_x^2 v+ \bar{f}(v)\partial_x v, \partial_x^2 v\>_H \\
 =& -\|\partial_x^2 v\|_{L^2}^2  -\int_0^1 \bar{f}(v)\partial_x  v \partial_x^2  v dx\\
\le & -\|\partial_x^2 v\|_{L^2}^2 +C\|\partial_x^2 v\|_{L^2} \|v\|_V\\
\le &-\frac{1}{2}\|\partial_x^2 v\|_{L^2}^2 + C\|v\|_V^2,  \  v\in H_n\subseteq V,
 \end{split}
\end{equation*}
i.e.  $(H3)$ also holds.

Therefore, the assertion follows from Theorem 2.2.
\end{proof}

\begin{rem}
(1) If we take $\bar{f}=0$ and $g(x)=\sum_{i=0}^{2n+1} a_i x^i$ with $a_{2n+1}<0$ (for some fixed $n\in \N$),
then (\ref{RDE}) reduces to the classical backward stochastic reaction-diffusion equations.

(2) If $g=0$, then (\ref{RDE})
covers the backward stochastic Burgers type equations.
\end{rem}
\subsection{Backward stochastic tamed 3D Navier-Stokes equation}

The last example  is  a tamed version of backward stochastic 3D Navier-Stokes equation. Stochastic tamed 3D Navier-Stokes equation
has been investigated  in a series of works of  R\"{o}ckner et al  \cite{RZ09a, RZ09, RZZ,RZ10}.
The classical 3D Navier-Stokes equations (i.e. $g_N=0, \  B=0$ in \ref{Tamed NSE}) is
a standard model to
 describe the evolution of velocity fields of an
incompressible fluid (cf.\cite{F08,Li96,T84}), the uniqueness and regularity of weak
 solutions are still  open problems  up to now.

The authors in \cite{QTY12,SY09} have studied the backward stochastic 2D Navier-Stokes equation.
To the best of our knowledge, there is no result about
 backward stochastic 3D Navier-Stokes equation,
the backward stochastic tamed 3D Navier-Stokes equation can be viewed as a regularized version of the classical
backward stochastic 3D Navier-Stokes equation and it can be formulated
as follows:
\begin{equation}\label{Tamed NSE}
\begin{split}
 &dX_t= -\left[ \nu \Delta X_t-(X_t\cdot \nabla)X_t+\nabla p(t) -g_N\left(|X_t|^2\right) X_t+h(t,X_t,Z_t)\right] dt
 + Z_t dW_t, \\
   &div (X_t)=0,   ~ \ X_T=\xi,\\
   &X_t|_{\partial \Lambda}=0,
\end{split}
\end{equation}
where $\nu>0$ is the viscosity constant,  $p$ is the (unknown)  pressure and   the taming function $g_N : \mathbb{R}_+\rightarrow
\mathbb{R}_+$ is smooth and satisfies for some $N > 0$,
$$\begin{cases}
  g_N(r)=0,& \text{if}\  r\le N,\\
 g_N(r)=(r-N)/\nu,&  \text{if}\  r\ge N+1,\\
 0\le g_N^\prime(r)\le C,&   r\ge 0.
\end{cases}
$$
 The main feature of  (\ref{Tamed NSE})
is that if there is a bounded smooth solution to the backward (stochastic) 3D Navier-Stokes equation,
then this smooth solution must also
satisfy this backward tamed equation
  for some  large enough $N$.

Let $\mathcal{P}$ be the orthogonal (Helmhotz-Leray) projection from
 $L^2(\L,\R^3)$ to $H^0$ (cf.\cite{T84,Li96,F08}).
 For any $u\in H^0$ and $v\in L^2(\L,\R^3)$ we have
$$ \<u, v\>_{H^0} :=  \<u, \mathcal{P} v\>_{H^0}=\<u, v\>_{L^2} . $$
We consider the following Gelfand triple:
$$         V:=   H^1 \subseteq     H:=  H^0\subseteq  V^*= (H^{1})^*,  $$
then it is well known that  the following operators
$$  A: W^{2,2}(\Lambda, \mathbb{R}^3)\cap V\rightarrow H, \ Au=\nu \mathcal{P} \Delta u ; $$
$$        F: \mathcal{D_F} \subset  H\times V\rightarrow H; \  F(u, v)=- \mathcal{P} \left[ (u\cdot \nabla) v   \right],  \ F(u):=F(u,u)  $$
can be extended to  the following  well defined operators:
$$     A: V\rightarrow V^*; \ F: V\times V\rightarrow V^*.              $$
Moreover,  we have
\begin{equation}\label{estimate of bilinear term}
     { }_{V^*}\<F(u,v),  w\>_V=  -   { }_{V^*}\<F(u,w),  v\>_V, \   { }_{V^*}\<F(u,v),  v\>_V=0, \ u,v,w\in V.
\end{equation}

Without loss of generality we may assume $\nu=1$. Now we show the existence and uniqueness of solutions
to (\ref{Tamed NSE}).
\begin{exa}
Suppose $\xi\in L^\infty_{\mathcal{F}_T}(\Omega;V)$ and $h:[0,T]\times V\times L_2(U,H)\times \Omega\rightarrow V^*$ satisfies (H0)-(H4), then
 (\ref{Tamed NSE}) has a unique adapted solution $(X,Z)\in L^\infty_\mathcal{F}(\Omega\times[0,T];V)\times L^2_\mathcal{F}(\Omega;L^2([0,T],L_2(U,H)))$.
Moreover, it satisfies that
$$\sup_{t\in[0,T]}\|X_t\|_V^2+\frac{1}{2}E\int_0^T\|Z_s\|_{L_2(U,V)}^2ds\leq C(1+\|f\|_{L_{\mathcal{F}}^\infty(\Omega;L^1([0,T]))}+\|\xi\|^2_{L_{\mathcal{F}_T}^\infty(\Omega;V)}), \quad a.s..$$
\end{exa}
\begin{proof}  It is well known that  (\ref{Tamed NSE}) can be  rewritten
 into the following variational form:
$$   dX_t=-\left[ AX_t+ F(X_t)- \mathcal{P}\left(g_N\left(|X_t|^2\right) X_t \right) +h(t,X_t,Z_t)\right] dt + Z_tdW_t, \ X_T=\xi. $$
It is easy to see that all eigenvectors
$\{ e_i, i=1,2,\cdots \}\subset H^2$ of $A$  constitute  an orthonormal basis of $H^0$ and an orthogonal set in $H^1$,
i.e. $(H0)$ holds.

By H\"{o}lder's inequality we  have  the following  estimate:
$$   \|\psi\|_{L^3(\Lambda;\mathbb{R}^3)}\le  \|\psi\|_{L^2(\Lambda;\mathbb{R}^3)}^{1/2}
\| \psi\|_{L^6(\Lambda;\mathbb{R}^3)}^{1/2}, \ \psi\in L^6(\Lambda; \mathbb{R}^3).   $$
Note that $W_0^{1,2}(\Lambda; \mathbb{R}^3) \subseteq L^6(\Lambda; \mathbb{R}^3)$,
 then by (\ref{estimate of bilinear term})  one can show  that
\begin{equation*}
 \begin{split}
 &~~ { }_{V^*}\<F(u)-F(v),u-v\>_V \\
&= -  { }_{V^*}\<F(u-v), v\>_V \\
&\le C  \|u-v\|_V \|u-v\|_{L^3(\Lambda;\mathbb{R}^3)}  \|v\|_{L^6(\Lambda;\mathbb{R}^3)} \\
 &\le C  \|u-v\|_V^{3/2} \|u-v\|_H^{1/2}  \|v\|_{L^6(\Lambda;\mathbb{R}^3)} \\
& \le \frac{1}{2} \|u-v\|_V^{2} + C \|v\|_{L^6(\Lambda;\mathbb{R}^3)}^4 \|u-v\|_H^{2},
\ u,v\in V.
 \end{split}
\end{equation*}
Hence we have the following estimate (recall that $\nu=1$):
$$  { }_{V^*}\<Au+F(u)-Av-F(v), u-v\>_V
 \le -\frac{1}{2} \|u-v\|_V^{2} + C\left(1+ \|v\|_{L^6(\Lambda;\mathbb{R}^3)}^4 \right) \|u-v\|_H^{2}. $$
By the definition of $g_N$ and (\ref{GN inequality}) we have
 \begin{equation*}
 \begin{split}
 & -{ }_{V^*} \<\mathcal{P}( g_N(|u|^2)u) - \mathcal{P}(g_N(|v|^2)v), u-v\>_V \\
=& - \< g_N(|v|^2) (u-v),  u-v\>_H + \< ( g_N(|v|^2)-  g_N(|u|^2) )  u,  u-v\>_H \\
\le &  \int_{\{ |u|>|v| \}  }   ( g_N(|v|^2)-  g_N(|u|^2) ) (|u|^2- u\cdot v) d x\\
 & +  \int_{\{ |u|\le |v| \}  }   ( g_N(|v|^2)-  g_N(|u|^2) ) (|u|^2- u\cdot v) d x\\
\le & C \int_{\{ |u|\le |v| \}  }   \left||v|^2- |u|^2 \right|\cdot|u|\cdot |u-v| dx\\
\le & C \int_{\{ |u|\le |v| \}  }   |u|^2\cdot |u-v|^2 dx\\
\le & C\|v\|^2_{L^6(\Lambda;\mathbb{R}^3)}  \|u-v\|^2_{L^3(\Lambda;\mathbb{R}^3)} \\
\le & C \|v\|^2_{L^6(\Lambda;\mathbb{R}^3)}  \|u-v\|_H  \|u-v\|_V \\
\le & \frac{1}{4} \|u-v\|_V^2 + C  \|v\|^4_{L^6(\Lambda;\mathbb{R}^3)} \|u-v\|^2_H , \ u,v\in V.
 \end{split}
\end{equation*}
Hence $(H2)$ holds with $\rho(v)=C  \|v\|_{L^6(\Lambda;\mathbb{R}^3)}^4$.

We recall the following estimate for $v\in span\{e_1, e_1, \cdots, e_n  \}$ (cf.\cite[Lemma 2.3]{RZ09a}):
\begin{equation}
\begin{split}
  \<Av, v\>_V&=\<\mathcal{P}\Delta v, (I-\Delta) v\>_{H}\le -\|v\|_{H^2}^2+ \|v\|_{V}^2 ;\\
\<F(v), v\>_V&=- \<\mathcal{P}\left[( v \cdot \nabla) v \right],  (I-\Delta)v\>_{H} \le \frac{1}{2} \|v\|_{H^2}^2+\frac{1}{2}\| |v|\cdot |\nabla v| \|_{H}^2 ;\\
-\<\mathcal{P}\left(g_N(|v|^2)v\right), v\>_V&=-\<\mathcal{P}\left( g_N(|v|^2)v \right), (I-\Delta)v\>_{H} \le
  -\| |v|\cdot |\nabla v| \|_{H}^2+ CN \|v\|_{V}^2.
\end{split}
\end{equation}
Then  it is  easy to  verify $(H3)$ as follows:
$$ \<Av+F(v)-\mathcal{P}(g_N(|v|^2)v),  v\>_V\le -\frac{1}{2}  \|v\|_{H^2}^2 + C(N+1) \|v\|_{V}^2,   \  v\in span\{e_1, e_1, \cdots, e_n  \}.   $$
Concerning the growth condition, we have   that
$$ \|F(v)\|_{V^*}\le C \|v\|_{L^4(\Lambda; \mathbb{R}^3)}^2\le C \|v\|_V^{2}, \  v\in V.$$
By \eref{GN inequality} we have
$$  \|g_N(|v|^2)v\|_{V^*}^2\le C \|v\|_{L^6(\Lambda;\mathbb{R}^3)}^2 \le C \|v\|_V^2,
 v\in V.   $$
Hence  we know that   $(H4)$.

Then  the   existence of a unique solution to (\ref{Tamed NSE}) follows from Theorem 2.2.
\end{proof}

\end{document}